\documentclass[11pt]{amsart}
\usepackage{amsmath,amsthm, amscd, amssymb, amsfonts}
\usepackage[all]{xy}
\usepackage{mathrsfs}
\usepackage{enumitem}

\usepackage[ansinew]{inputenc}
\usepackage{graphicx,fancyhdr}

\usepackage[dvips, dvipsnames, usenames]{color}

\newcommand{\comment}[1]{}

\numberwithin{equation}{section}
\newtheorem{theorem}{Theorem}[section]
\newtheorem{lemma}[theorem]{Lemma}

\newtheorem{prop}[theorem]{Proposition}

\theoremstyle{definition}
\newtheorem{definition}[theorem]{Definition}
\newtheorem{example}[theorem]{Example}

\theoremstyle{remark}
\newtheorem{remark}[theorem]{Remark}

\def\pf{\begin{proof}}
\def\epf{\end{proof}}

\newcommand\GK{\operatorname{GK-dim}}
\newcommand\vect{\operatorname{vect}}

\newcommand{\ku}{ \Bbbk}
\newcommand{\kut}{ \Bbbk^{\times}}

\newcommand\I{\mathbb I}

\newcommand\N{\mathbb N}

\newcommand{\bq}{\mathbf{q}}

\newcommand\Tb{\mathbb T}

\newcommand\Z{\mathbb Z}
\def\zt{\Z^{\theta}}

\renewcommand{\_}[1]{_{\left( #1 \right)}}
\renewcommand{\^}[1]{^{\left( #1 \right)}}
\newcommand{\bj}[1]{_{[ #1 ]}}

\newcommand\cB{\mathcal{B}}

\def\dpn{\widetilde{\mathcal{B}}}
\newcommand{\dpnq}{\dpn(V)}

\newcommand{\lu}{\mathcal{L}}
\newcommand{\luq}{\lu(V)}

\newcommand\E{\mathcal{E}}

\newcommand\cI{\mathcal{I}}
\newcommand\cJ{\mathcal{J}}

\newcommand\cM{\mathcal{M}}

\newcommand\cP{\mathcal{P}}

\newcommand{\Ss}{{\mathcal S}}
\newcommand{\Sso}{\overline{\mathcal S}}
\newcommand\T{\mathcal{T}}
\newcommand\Uc{\mathcal{U}}
\newcommand\ug{\mathfrak{u}}
\newcommand\uq{\mathfrak{u}(V)}

\newcommand\bd{\mathfrak{bd}}

\newcommand\g{\mathfrak g}
\newcommand\gl{\mathfrak{gl}}

\newcommand\tg{\mathfrak t}

\newcommand\ad{\operatorname{ad}}
\newcommand{\Alg}{\Hom_{\text{alg}}}
\newcommand\Aut{\operatorname{Aut}}

\newcommand\car{\operatorname{char}}
\newcommand\Bider{\operatorname{Bider}}
\newcommand\Der{\operatorname{Der}}
\newcommand\End{\operatorname{End}}

\newcommand\id{\operatorname{id}}

\newcommand\gr{\operatorname{gr}}

\newcommand{\Hom}{\operatorname{Hom}}

\def\ydh{{}^{H}_{H}\mathcal{YD}}
\def\ydhb{{}^{H^\text{bop}}_{H^\text{bop}}\mathcal{YD}}

\def\hyd{{\mathcal{YD}}^{H}_{H}}

\def\kyd{{\mathcal{YD}}^{K}_{K}}

\def\ydk{{}^{K}_{K}\mathcal{YD}}
\def\ydl{{}^{L}_{L}\mathcal{YD}}

\newcommand{\toba}{\mathcal{B}}

\newcommand{\ot}{\otimes}

\newcommand{\ydG}{{}^{\ku \Gamma }_{\ku \Gamma }\mathcal{YD}}

\newcounter{tabla}\stepcounter{tabla}

\begin{document}


\title[unrolled hopf algebras]{On unrolled Hopf algebras}

\author[Andruskiewitsch;  Schweigert]
{Nicol\'as Andruskiewitsch and Christoph Schweigert}

\address{FaMAF-CIEM (CONICET), Universidad Nacional de C\'ordoba,
Medina Allende s/n, Ciudad Universitaria, 5000 C\' ordoba, Rep\'
ublica Argentina.} \email{andrus@famaf.unc.edu.ar}

\address{Fachbereich Mathematik, Universit\"at Hamburg,
Bereich Algebra und Zahlentheorie\\
Bundesstra\ss e 55, \ D\,--\,20\,146\, Hamburg
} \email{christoph.schweigert@uni-hamburg.de}

\thanks{The work of N. A. was partially supported by CONICET, Secyt (UNC) and the
MathAmSud project GR2HOPF. 
C. S. was partially supported by the Collaborative Research Centre 676 ``Particles,
Strings and the Early Universe - the Structure of Matter and Space-Time'', and by the RTG 1670
``Mathematics inspired by String theory and 
Quantum Field Theory''.
Most of this work was done during a visit of N. A. to the
University of Hamburg, supported by the Alexander von
Humboldt Foundation, in April 2015.
}

\begin{abstract} 
We show that the definition of unrolled Hopf algebras can be
naturally extended to the Nichols algebra $\toba(V)$ of
a Yetter-Drinfeld module $V$ on which a Lie algebra $\g$ 
acts by biderivations. Specializing to Nichols algebras of 
diagonal type, we find unrolled versions of the 
small, the De Concini-Procesi and the Lusztig divided power quantum group,
respectively.
\end{abstract}

\maketitle

\section{Introduction}\label{sect:intro}
\subsection{}\label{subsect:main}
In the recent papers \cite{CGP, GPT}, a so called unrolled version of quantum $sl(2)$ was introduced,  with applications to quantum topology; the definition was generalized to simple finite-dimensional Lie algebras in \cite{GP}. In the present article, we propose a generalization of this notion and
embed it into the appropriate conceptual context.

Recall that the  unrolled quantum $sl(2)$ is defined 
as the smash product of $U_q(sl(2))$ by the universal 
enveloping algebra of the Lie algebra of dimension 1.
Our starting point is the observation in Lemma \ref{lemma:bider}: 
given an action of the universal enveloping algebra 
$U(\g)$ of a Lie algebra $\g$ on a Hopf algebra $H$, 
the smash product is a Hopf algebra, 
if and only if $\g$ acts on $H$ by biderivations.
We next observe that,  if $V$ is a Yetter-Drinfeld module 
over a group $G$, then the Lie algebra
$\bd_V:=\End_G^G(V)$ of endomorphisms of the Yetter-Drinfeld 
module $V$ acts by biderivations on the Nichols algebra 
$\toba(V)$. Hence,  we can form the Hopf algebra
$(\toba(V)\#\ku G) \rtimes U(\bd_V)$ which we call the \emph{ unrolled bosonization} of $V$.
If $\dim V$ is finite, then its Gelfand-Kirillov dimension
can be expressed in terms of the Gelfand-Kirillov dimension
of $\toba(V)$ and the dimension of $\bd_V$.

The construction of unrolled bosonizations 
extends to a Lie subalgebra $\g$ of $\bd_V$, pre- or post-Nichols algebras (in the place of $\toba(V)$), and to deformations thereof, provided
that the action of the Lie algebra $\g$ preserves the relevant defining relations. In particular, we  define the unrolled version of the 
quantum double of a finite-dimensional Nichols algebra of diagonal type. 

\subsection{Preliminaries}\label{subsect:preliminaries}
Fix a field $\ku$ and let $H$ be a Hopf algebra over $\ku$. We use standard notation: $\Delta$, $\varepsilon$, $\Ss$, $\Sso$ are respectively the comultiplication, the counit, the antipode (always assumed to be bijective) and
the inverse of the antipode.

We denote by $\ydh$ the category of Yetter-Drinfeld modules over $H$ as in \cite{AS Pointed HA}. For $V, W \in \ydh$, we denote by $\Hom_H^H (V, W)$, $\End_H^H (V)$, $\Aut_H^H (V)$
the spaces of morphisms, respectively endomorphisms, automorphisms in $\ydh$.
Let $R$ be a Hopf algebra in the braided monoidal category $\ydh$,
with comultiplication denoted by
$r \mapsto r\^1 \ot r\^2$.
Recall that the \emph{bosonization} $R \# H$ is the Hopf algebra over
$\ku$  with
underlying vector space $R\otimes H$, smash product multiplication and smash coproduct
comultiplication; i.e., for all $r,s \in R$, $a,b \in H,$
\begin{align}\label{eq:smash product}
(r \# a)(s \# b) &= r (a\_{1} \cdot s) \# a\_2b,
\\\label{eq:smash coproduct}
\Delta(r \# a) &= r\^1 \# (r\^2){\_{-1}} a\_1 \ot (r\^2){\_0} \# a\_2.
\end{align}
Here we write $r\# h$ for $r\ot h$.

We also introduce the category $\hyd = \ydhb$ 
of \emph{right-right} Yetter-Drinfeld modules over $H$.
Thus $M \in \hyd$ means that $M$ is a right $H$-module and a right $H$-comodule (with coaction $\varrho$),
and satisfies the compatibility axiom
\begin{align}\label{eq:yd-right-compatible}
\varrho(m \cdot h) &= m\_{0} \cdot h\_2 \ot \Ss(h\_1) m\_1 h\_3, & m &\in M, \, h\in H.
\end{align}
The tensor category $\hyd$ is braided, with braiding
$c(m \ot n) = n\cdot m\_1 \ot m\_0$, for all $m \in M$, $n\in N$, $M, N\in \hyd$.
For right-right Yetter-Drinfeld modules $V, W \in \hyd$, 
we use the notion $\Hom_H^H (V, W)$, $\End_H^H (V)$, $\Aut_H^H (V)$ are
as before.

Let $T$ be a Hopf algebra in the braided monoidal category $\hyd$ of
right-right Yetter-Drinfeld modules,
with comultiplication denoted by
$t \mapsto t\^1 \ot t\^2$.
In this case, the \emph{bosonization} $H \# T$ is the Hopf algebra over
$\ku$  with
underlying vector space $H\otimes T$, smash product multiplication and smash coproduct
comultiplication; i.e.
\begin{align}\label{eq:smash product-right}
(a \# t)(b \# u) &=  a b\_{1} \#  (t \cdot b\_2) u,
\\\label{eq:smash coproduct-right}
\Delta(a \# t) &= a\_1 \# (t\^1){\_{0}}  \ot a\_2 (t\^1){\_{1}} \# t\^2,
\end{align}
for all $t,u \in R$, $a,b \in H$. Here we write $h\# t$ for $h\ot t$.

If $\Gamma$ is an abelian group, then we denote by
$\ku_g^{\chi}$ the one-dimensional object in $\ydG$ with coaction given
by the group element $g\in\Gamma$ and action given by the character
$\chi \in \widehat\Gamma$. For a Yetter-Drinfeld module
$V\in \ydG$, the corresponding isotypic component is denoted by
$V_g^{\chi}$. A Yetter-Drinfeld module has a natural structure
of a braided
vector space. For a braided vector space $V$, denote by
$\toba(V)$ its Nichols algebra and by $\cJ = \cJ(V)$ its 
ideal of defining relations, cf. \cite{AS Pointed HA}; 
so that $\toba(V) \simeq T(V)/\cJ(V)$.

\subsection*{Acknowledgements} N.A. thanks I. Angiono for some interesting exchanges.

\section{Unrolled Hopf algebras}\label{sect:lifting-diagonal-Ug}
\subsection{}\label{subsubsect:smash}

Let $L$ be a Hopf algebra. Recall that a (left) $L$-module algebra is an algebra $A$ which is also
an  $L$-module with action  $\cdot:L \otimes A \to A$ such that for all $\ell \in L$ and all  $a,b  \in A$ the compatibility conditions 
\begin{align}\label{eq:smash1}
\ell \cdot (ab) &= (\ell\_1 \cdot a)(\ell\_2 \cdot b),
\\ \label{eq:smash-unit}
\ell \cdot 1 &= \varepsilon(\ell)1.
\end{align}
for product and unit hold.
It is well-known that \eqref{eq:smash1} and \eqref{eq:smash-unit}  mean that $A$ is an algebra in the monoidal category
$_L\cM$ of left $L$-modules.

In this paper, we are interested in the case of a Hopf algebra $H$
that is also an $L$-module algebra, where $L$ is a Hopf algebra as well.
In this case,  we impose the following consistency conditions:
\begin{align}
\label{eq:smash2}
\Delta(\ell \cdot a) &= \ell\_1 \cdot a\_1\otimes \ell\_2 \cdot a\_2,
\\\label{eq:smash-ep}
\varepsilon(\ell \cdot a) &= \varepsilon(\ell) \varepsilon(a),
\\ \label{eq:smash3}
\ell\_1  \otimes\ell\_2 \cdot a &=  \ell\_2  \otimes\ell\_1 \cdot a,
\end{align}
for all $\ell \in L$ and all $a,b  \in H$.
Then $H\rtimes L := H\otimes L$ with the tensor product structure
as a coalgebra and with the smash product
\eqref{eq:smash product} for the algebra structure is a Hopf algebra;
see \cite{M}, \cite[1.2.10]{AN} (in this second paper a different notation is used).
We shall say that $H$ is a \emph{$L$-module Hopf algebra}.

\begin{remark}
The following perspective shows that it is natural to impose
these consistency conditions.
The category $_L\cM$ of left $L$-modules is monoidal, but not
braided; thus $H$ cannot be interpreted as a Hopf algebra in 
$_L\cM$.
Still, it can be interpreted in terms of monads. Recall that
$A$ has the structure of an algebra in the monoidal category
$_L\cM$ of left $L$-modules, if and only if the endofunctor $T: {}_L\cM \to {}_L\cM$,
$T(X) = A \otimes X$ has the structure of a monad.

Also recall \cite{BLV} that
a bimonad structure on a monad $T$ on a monoidal category
consists of a comonoidal structure on the functor $T$, i.e.\
a natural transformation
$$T_2: T(X \ot Y) = H \ot (X \ot Y) \to T(X) \ot T(Y) = (H\ot X)
\ot (H \ot  Y),$$
and a morphism $T_0: T(1)  \to1$. They have to obey axioms
generalizing coassociativity and counitality. If
$H$ is a bialgebra in a braided monoidal category, the monad
$T(-)=H\otimes -$ can be endowed
via the coproduct $\Delta: H \to H \ot H$ with the natural
transformation
$$ T_2(a \ot x \ot y) =  (a\_1 \ot x) \ot (a\_2 \ot  y)\,\,, $$
where we used Sweedler notation for $\Delta$. The morphism
$T_0$ is induced from the  counit $\varepsilon: H \to\ku$.

Now let $L$ be another Hopf algebra and $H$ be an $L$-module
algebra.
The fact that $T_2$ is a morphism in $_L\cM$ is then equivalent
to the consistency conditions \eqref{eq:smash2} and \eqref{eq:smash3},
while condition \eqref{eq:smash-ep} amounts to the fact
that $\varepsilon$ is a morphism in $_L\cM$.
Thus $T(-)=H\otimes -$ is a bimonad on the monoidal category
$_L\cM$, if and only if
the requirements \eqref{eq:smash2}, \eqref{eq:smash-ep},
and \eqref{eq:smash3} hold. It is a Hopf monad, if and only if
$H$ is a Hopf algebra.
The  Hopf monad in Vec$_\ku$ (i. e., Hopf algebra) $H\rtimes L$
corresponds to the forgetful functor as described in \cite[Proposition 4.3]{BLV}.
\end{remark}

\begin{remark}  Here is another way to interpret $H\rtimes L$,  dual to \cite[1.1.5]{AN}.
Let $H$ be a $L$-module Hopf algebra. Then $H$, endowed with the trivial coaction, is a Hopf algebra in $\ydl$
and $H\rtimes L \simeq H\# L$. Indeed, \eqref{eq:smash3} is equivalent to the compatibility in $\ydl$.
\end{remark}

\subsection{}\label{subsec:duals}
Now turn to the situation of two Hopf algebras
$H$ and $U$, provided with a non-degenerate bilinear form
$(\, \vert \,): H \ot U \to \ku$. We extend this bilinear form
to  a non-degenerate bilinear form $(\, \vert \,): H \ot H \ot U \ot U \to \ku$ by
\begin{align}\label{eq:duality}
(a\ot \widetilde a \vert u\ot \widetilde u) &:= (a\vert \widetilde{u})(\widetilde{a}\vert u), &\text{ for }\,
a, \widetilde a \in H,\,  u, \widetilde u \in U.
\end{align}
We assume that the pairing $(\, \vert \,)$ is such that for every $a, \widetilde a \in H$, $u, \widetilde u \in U$, the following identities hold
\begin{align} \label{eq:dual1}
(a \widetilde a \vert u) &= (a \ot  \widetilde a \vert \Delta(u)) = (a\vert u\_2)(\widetilde{a}\vert u\_1), & (1\vert u) &= \epsilon (u), \\ \label{eq:dual2}
(a  \vert u\widetilde u) &= (\Delta(a)  \vert  u \ot\widetilde u ) = (a\_2\vert u)(a\_1\vert\widetilde{u}), & (a\vert 1) &= \epsilon (a), \\ \label{eq:dual3}
(\Ss(a)  \vert u) &= (a \vert \Ss(u)). &&
\end{align}
Such a pairing is called a Hopf pairing on $H$ and $U$.

\begin{lemma}\label{pairinglemma}
Assume that the two Hopf algebras $H$ and $U$ are $L$-modules 
and that there is
a Hopf pairing on $H$ and $U$. Assume that the pairing is compatible
with the $L$-action involving the antipode of $L$,
        \begin{align} \label{eq:duality-modL}
(\ell \cdot a  \vert u) &= (a \vert \Ss(\ell) \cdot u), & a  \in H,\,  u \in U,\,  \ell  \in L.
        \end{align}
Then the Hopf algebra $H$ is an $L$-module Hopf algebra, 
if and only if $U$ is so.
\end{lemma}

\pf
Let $\ell \in L$, $u,v  \in U$ and $a \in H$.  We compute
\begin{align*}
( a\vert \ell \cdot (uv)) &= ( \Sso(\ell) \cdot a\vert uv) =  ((\Sso(\ell) \cdot a)\_2\vert u)((\Sso(\ell) \cdot a)\_1\vert v);
\\
( a\vert(\ell\_1 \cdot u)(\ell\_2 \cdot v)) &= ( a\_2 \vert \ell\_1 \cdot u)( a\_1\vert \ell\_2 \cdot v)
= (\Sso(\ell\_1) \cdot a\_2 \vert u)(\Sso(\ell\_2) \cdot a\_1\vert v)
\\
&= (\Sso(\ell)\_2 \cdot a\_2 \vert u)(\Sso(\ell)\_1 \cdot a\_1\vert v).
\end{align*}
Hence \eqref{eq:smash1} holds for $U$  if and only if
$( a\vert \ell \cdot (uv)) = ( a\vert(\ell\_1 \cdot u)(\ell\_2 \cdot v))$ for all $\ell \in L$, $u,v  \in U$, $a \in H$,  if and only if
$((\widetilde{\ell} \cdot a)\_2\vert u)((\widetilde{\ell} \cdot a)\_1\vert v) = (\widetilde{\ell}\_2 \cdot a\_2 \vert u)(\widetilde{\ell}\_1 \cdot a\_1\vert v)$for all $\widetilde{\ell} \in L$, $u,v  \in U$, $a \in H$, if and only if \eqref{eq:smash2} holds for $H$. Thus
\eqref{eq:smash1} holds for $H$ if and only if \eqref{eq:smash2} holds for $U$.

Similarly \eqref{eq:smash-unit} holds for $U$ if and only if \eqref{eq:smash-ep} holds for $H$ and vice versa.
Finally, \eqref{eq:smash3} holds for $H$ if and only if it holds for $U$:
\begin{gather*}
\ell\_1  \otimes\ell\_2 \cdot u =  \ell\_2  \otimes\ell\_1 \cdot u, \quad \forall u \iff \Sso(\ell\_1) (a \vert \otimes\ell\_2 \cdot u)
=  \Sso(\ell\_2)  (a \vert \ell\_1 \cdot u), \quad \forall u, a \iff
\\   \Sso(\ell\_1) (\Sso(\ell\_2) \cdot a \vert  u)
 =  \Sso(\ell\_2)  ((\Sso(\ell\_1) \cdot a \vert u), \quad \forall u, a
\iff \Sso(\ell)\_2 (\Sso(\ell)\_1 \cdot a \vert  u)
\\ =  \Sso(\ell)\_1  ((\Sso(\ell)\_2 \cdot a \vert u,) \quad \forall u, a
\iff \Sso(\ell)\_2 \ot \Sso(\ell)\_1 \cdot a
=  \Sso(\ell)\_1  \ot \Sso(\ell)\_2 \cdot a,  \quad \forall a.
\end{gather*}
\epf

\subsection{}\label{subsec:smash-braided}
We next extend our construction to Hopf algebras in braided
monoidal categories. To this end,
let now $K$ be a Hopf algebra, $\toba$ a Hopf algebra  in
the braided category $\ydk$. Let $L$\ be another Hopf algebra as before,
and assume that $\toba$ is also an $L$-module algebra.
We extend the action of the Hopf algebra $L$ to
the bosonization $H := \toba \# K$ by
$\ell \cdot (b \# k)  := (\ell \cdot b) \# k$, for $\ell \in L$,
$b \in \toba$ and $k \in K$. 

Then straightforward verifications show that

\begin{itemize} [leftmargin=*]
\item 
The bosonization $H$ is a $L$-module algebra $\iff$ 
The actions of $L$ and $K$ on $\toba$ commute.

\item \eqref{eq:smash-ep} holds for $H \iff$ \eqref{eq:smash-ep} holds 
for $\toba$.\\
From now on, we assume that this is the case.

 \item \eqref{eq:smash2} holds for $H \iff$ \eqref{eq:smash2} holds for $\toba$ and
 the action of $\ell$ on $\toba$ is a morphism of $K$-comodules for all $\ell\in L$.

  \item \eqref{eq:smash3} holds for $H \iff$ \eqref{eq:smash3} holds for $\toba$.
\end{itemize}

In other words,  the action of $L$ on the bosonization 
$H = \toba \# K$ satisfies \eqref{eq:smash-ep}, \eqref{eq:smash2}
and \eqref{eq:smash3}, if and only if so does the action of $L$ on $\toba$, and the homothety $\eta_\ell$ for $\ell\in L$ is a morphism of
Yetter-Drinfeld modules,
$\eta_\ell \in \End^K_K \toba$ for all $\ell\in L$. This leads to

\begin{definition}
An  \emph{$L$-module braided Hopf algebra} is a  Hopf algebra $\toba$ in
the braided category $\ydk$ that is also  a $L$-module algebra, that satisfies \eqref{eq:smash-ep}, \eqref{eq:smash2} and  \eqref{eq:smash3},  and such that the homothety
$\eta_\ell \in \End^K_K \toba$ for all $\ell\in L$. 
\end{definition}

We have just seen:
for an $L$-module braided Hopf algebra, the 
bosonization $H := \toba \# K$ is an $L$-module Hopf algebra
over $\ku$ and we can form the Hopf algebra $H \rtimes L = (\toba \# K)  \rtimes L$.

\medbreak
As in subsection \ref{subsec:duals}, we consider the situation with
non-degenerate pairings; this time internal to the braided monoidal
category $\ydk$ instead of $\vect \ku$. Concretely,
let $\E$ be another Hopf algebra  in the category $\ydk$ provided with a non-degenerate bilinear form  $(\, \vert \,): \toba \ot \E \to \ku$, and  extend it  by \eqref{eq:duality}
to a pairing $\toba\otimes\toba\otimes\E\otimes \E\to\ku$.
\begin{itemize} [leftmargin=*]\renewcommand{\labelitemi}{$\diamond$}
\item  
The fact that the pairing is internal to the category 
$\ydk$ means that 
the bilinear form $(\, \vert \,)$ is a morphism in
the monoidal category $\ydk$, where $\ku$ is endowed with 
the structure of a trivial Yetter-Drinfeld module.

\item 
We assume that
for every $a, \widetilde a \in \toba$, $u, \widetilde u \in \E$,
the conditions \eqref{eq:dual1}, \eqref{eq:dual2} and \eqref{eq:dual3}
of a Hopf pairing, relating coproduct, product, unit and 
counit of $\toba$ and $\E$ hold.
\end{itemize}

Then we have in the braided category $\ydk$ exactly the same situation
we considered in lemma \ref{pairinglemma} in the braided category
$\vect \ku$. The same calculations, this time in the
category $\ydk$, yield:

\begin{lemma}\label{pairinglemma:braided}
Assume that both $\toba$ and $\E$ are $L$-modules and that 
condition \eqref{eq:duality-modL} on the Hopf pairing
$(\, \vert \,)$ holds.
Then $\toba$ is a $L$-module braided Hopf algebra, if and only if $\E$ is so. \qed
\end{lemma}

\subsection{}\label{subsec:smash-Liealg}
Let $\g$ be a Lie algebra over the field $\ku$. We specialize
to $L$-module braided Hopf algebras where the Hopf algebra $L = U(\g)$ is
the universal enveloping algebra of $\g$.
Then the conditions \eqref{eq:smash1} and \eqref{eq:smash-ep}
in the definition of an $L$-module Hopf algebra $H$
just mean that $\g$ acts on $H$ by $\ku$-derivations, while
condition \eqref{eq:smash3} is for free, 
due to the cocommutativity of $U(\g)$. 
Condition \eqref{eq:smash2} amounts to the condition
\begin{align}\label{eq:smash2-lie}
\Delta(x\cdot a) &= x \cdot a\_1\otimes  a\_2 + a\_1\otimes x \cdot a\_2,&
\varepsilon(x \cdot a) &= 0, & &
\end{align}
for all $x\in \g$ and $a\in H$. In other words, condition \eqref{eq:smash2-lie}
tells us that
$\g$ acts on $H$ by $\ku$-coderivations. We summarize all
conditions by saying that $\g$ acts on $H$ by \emph{$\ku$-biderivations}:
$\g$ acts by endomorphisms that are simultaneously  
$\ku$-derivations and $\ku$-coderivations.
Thus we have:

\begin{lemma}\label{lemma:bider}
Let $H$ be a Hopf algebra and let $\g$ be a Lie algebra acting on $H$ by $\ku$-biderivations.
Then $H$ is a $U(\g)$-module Hopf algebra and we can form 
the Hopf algebra $H \rtimes U(\g)$. \qed
\end{lemma}

The following remarks on biderivations are useful:
\begin{itemize} [leftmargin=*]\renewcommand{\labelitemi}{$\diamond$}
\item 
For any Hopf algebra $H$, the subspace
$\Bider_{\ku} (H) := \{x\in \Der_{\ku} (H): x \text{ is a coderivation}\}$ is a Lie subalgebra of $\Der_{\ku}(H)$.
\item
If $x\in \Der (H)$ and if $a, b \in H$ fulfill \eqref{eq:smash2-lie} for $x$, then so does their product $ab$. Hence it is enough to check the biderivation
property \eqref{eq:smash2-lie} for a given derivation $x$ on a family
of generators of $H$.
\end{itemize}

\begin{remark} Let $H$ be a Hopf algebra and let $\g$ be a Lie algebra acting on $H$ by $\ku$-coderivations.
Let $H_0$ be the coradical, and $(H_n)_{n\ge 0}$ the coradical filtration, of $H$.
If $H_0$ is $\g$-stable, then $H_n$ is $\g$-stable for all $n\ge 0$ by the defining condition \eqref{eq:smash2-lie}.
Hence $\g$  acts on $\gr H$ by $\ku$-coderivations.

Assume that $H_0$ is a Hopf subalgebra, that $\g$  acts on $H$ by $\ku$-biderivations and that $H_0$ is $\g$-stable.
Then $\g$  acts on  the graded object $\gr H$ by $\ku$-biderivations.

Notice that $\g$  may act on $H$ by $\ku$-biderivations with $H_0$ not being $\g$-stable.
For instance, let $x \in H$ primitive. Then $D = \ad x$ is a $\ku$-biderivation.
If there exists $g \in G(H)$ such that $gx = qxg$ with $q\in \kut - \{1\}$, then $D(g) = (1-q) xg \notin H_0$.

\end{remark}

\subsection{}\label{subsec:pointed-smash}
In this context, suppose that $H$ is pointed
and set $G := G(H)$
the group of group-like elements of $H$.
Let $\g$ act on $H$ by derivations;  assume that $\g$ acts  trivially on $\ku G$.
Let $g, t \in G$ and $\cP_{g,t}(H)  := \{a\in H: \Delta(a) = g\otimes a + a \otimes t\}$
the space of $(g,t)$ skew-primitive elements.
Then the coderivation property  \eqref{eq:smash2-lie} implies that
$ \cP_{g,t}(H)$ is a $\g$-submodule for all  $g, t\in G$.
Summarizing, we have

\begin{lemma}
Let $\g$ be a Lie algebra acting by derivations on a pointed Hopf algebra $H$, $G = G(H)$.
Assume that
\begin{itemize}
 \item $\g$ acts trivially on $\ku G$.
 \item $H$ is generated by group-like and skew-primitive elements.
\end{itemize}

Then the following are equivalent:
\begin{enumerate}\renewcommand{\theenumi}{\alph{enumi}}
 \item $\g$ acts on $H$ by $\ku$-biderivations, i.e. \eqref{eq:smash2-lie} holds.
 \item $ \cP_{g,t}(H)$ is a $\g$-submodule for all  $g, t\in G$.
 \item $ \cP_{g,1}(H)$ is a $\g$-submodule for all  $g\in G$. \qed
\end{enumerate}
\end{lemma}

\subsection{}\label{subsec:smash-nichols-general} Let $K$ be a Hopf algebra and $V \in \ydk$. It is well-known that every
$d\in \Hom(V, T(V))$ extends uniquely to a derivation $D\in \Der(T(V))$
on the tensor algebra $T(V)$ by $D(1) = 0$
and
\begin{align}\label{eq:der-TV}
D_{\vert T^n(V)} &= \sum_{1\le j \le n} \id_{T^{j-1}(V)}\ot d \ot \id_{T^{n-j}(V)},
\end{align}
for $n > 0$. Thus every  Lie algebra map $\g \to \End(V)$ extends to a Lie algebra map $\g \to \Der(T(V))$.

\begin{prop}\label{lemma:unrolled-bosonization}
Let $V\in\ydk$.
Every  morphism of Lie algebras $\g \to \End_K^K(V)$ extends to
an action of the universal enveloping algebra $U(\g)$ on 
$T(V)\#K$ and to an action on $\toba(V)\#K$,
giving rise to the Hopf algebras $(T(V)\#K) \rtimes U(\g)$ and $(\toba(V)\#K) \rtimes U(\g)$.
\end{prop}

\pf
As explained, the action of $\g$ on $V$ extends uniquely
to an action of $\g$  on the tensor algebra $T(V)$ by derivations. Formula
\eqref{eq:der-TV} and the assumptions imply that this action
is by morphisms in the category $\ydk$.
By definition, \eqref{eq:smash2} holds in $V$, hence it holds in $T(V)$.
By \S \ref{subsec:smash-braided}, the action extended to $T(V)\#K$ satisfies the requirements
in \S \ref{subsubsect:smash}, hence we can form $(T(V)\#K) \rtimes U(\g)$. Second,
the action of $\g$ on $T^n(V)$ commutes with that of the braid group $\mathbb B_n$; since
the kernel of the  projection $T^n(V)\to \toba^n (V)$ is the kernel of the quantum symmetrizer,
$\g$ acts on the Nichols algebra $\toba(V)$ with the 
desired requirements.
\epf

\begin{definition}
Let $K$ be a Hopf algebra, $V \in \ydk$ and $\g$ a Lie subalgebra of $\bd_V:=\End_K^K(V)$. We call the Hopf algebra $(\toba(V)\#K) \rtimes U(\g)$
the\emph{ unrolled bosonization} of the Nichols algebra of $V$ by $\g$.
\end{definition}
One may define unrolled versions of bosonizations of pre-Nichols or post-Nichols algebras, see e.g \cite{AAR},
or of deformations of Nichols algebras, provided that the ideals of defining relations are preserved
by the action of $\bd_V$, or if $\bd_V$ is replaced by a suitable subalgebra.

\subsection{Finite $\GK$}

Our main reference for this subsection is \cite{KL}. Let $A$ be an associative $\ku$-algebra.
We say that a finite-dimensional subspace $V\subseteq A$ is
\textit{GK-deterministic} if
\begin{align*}
\GK A=\lim_{n \to \infty} \log_n \dim \sum_{0\le j\le n}V^n.
\end{align*}

\begin{lemma}\label{lemma:GKdim-smashproduct} \cite[Lemma 2.2]{AAH}
 Let $K$ be a Hopf algebra, $R$ a
 Hopf algebra in $\ydk$, $A$ a $K$-module algebra and
 $B$ an $R$-module algebra in $\ydk$. Assume that the actions of $K$ on $A$, of
 $K$ on $B$, of $K$ on $R$, and of $R$ on $B$ are locally finite.

 \begin{enumerate}[leftmargin=*,label=\rm{(\alph*)}]
        \item\label{item:GKdim-smashproduct}
        $\GK A\#K \le \GK A+\GK K$. If either $K$ or $A$ has a GK-deterministic subspace, then $\GK A\# K = \GK A + \GK K$.
        
        \medbreak
        \item\label{item:GKdim-smashproduct-braided}
        $\GK B\# R \le \GK B + \GK R$.
        If either $R$ or $B$ has a GK-deterministic subspace, then
        $\GK B\# R = \GK B + \GK R$.    \qed
 \end{enumerate}

\end{lemma}

Clearly,  a finite-dimensional Lie algebra $\g$ is a GK-deterministic subspace of $U(\g)$.
Thus we have:

\begin{example}
        Let $H$ be a Hopf algebra and let $\g$ be a Lie subalgebra of $\Bider_{\ku} (H)$ such that  $\GK H$, $\dim \g < \infty$. If
        the action of $\g$ on $H$ is locally finite, then
        \begin{align}
        \GK (H\rtimes U(\g)) = \GK H + \dim \g < \infty.
        \end{align}
\end{example}

Here are some particular cases:

\begin{itemize}[leftmargin=*]\renewcommand{\labelitemi}{$\circ$}
        \item If $H$ is a finite-dimensional Hopf algebra and  $\g$ is a Lie subalgebra of $\Bider_{\ku} (H)$,
        then $$\GK (H\rtimes U(\g)) =  \dim \g < \infty.$$
        
        \item Let $K$ be a Hopf algebra, $V\in \ydk$, $\g$ a Lie subalgebra of $\bd_V$, $\toba \in \ydk$ a pre-Nichols algebra of $V$ and $\E \in \ydk$ a post-Nichols algebra of $V$. Assume that the action of $\g$ descends to $\toba$  and $\E$, 
        \begin{align*}
        \GK K &< \infty,& \dim V &< \infty,& \GK \toba &< \infty,& \GK \E &< \infty.
        \end{align*}
        Clearly, $\dim \g < \infty$ and $\g$ acts locally finitely on $\toba \# K$ and $\E \# K$.
        If either $K$ or $\toba$, respectively $\E$, have a GK-deterministic subspace, then
\begin{align*}
        \GK \left((\toba \# K)  \rtimes U(\g)\right) = \GK \toba + \GK K + \dim \g < \infty, \\
        \GK \left((\E \# K)  \rtimes U(\g)\right) = \GK \E + \GK K + \dim \g < \infty.
\end{align*}
\end{itemize}

\section{The dual construction}\label{sect:dual-unrolled}
\subsection{}\label{subsubsect:cosmash}

Let $J$ be a Hopf algebra. A  $J$-comodule coalgebra is a coalgebra $C$ which is also
a right  $J$-comodule with coaction  $\varrho:C \to C \otimes J$, $\varrho (c) = c_{[0]} \ot c_{[1]}$,
 and counit $\varepsilon_C$ such that 
for all $c \in C$
\begin{align}\label{eq:smash-coprod}
(c\_1)\bj{0} \ot (c\_2)\bj{0} \ot (c\_1)\bj1 (c\_2)\bj1 &= (c\bj{0})\_1 \ot (c\bj{0})\_2 \ot c\bj1,
\\ \label{eq:smash-counit}
\varepsilon_C(c\bj{0}) c\bj1 &= \varepsilon_C(c).
\end{align}
Here \eqref{eq:smash-coprod} and \eqref{eq:smash-counit}  mean that $C$ is a coalgebra in the monoidal category
$\cM^J$ of right $J$-comodules. Assume that $C = H$ is a Hopf algebra and a  $J$-comodule coalgebra that satisfies:
\begin{align}
\label{eq:smash2-coprod}
(ab)\bj{0} \ot (ab)\bj1 &= a\bj{0} b\bj{0} \ot a\bj1b\bj1,
\\\label{eq:smash-ep-coprod}
\varrho(1) &= 1 \ot 1,
\\ \label{eq:smash3-coprod}
a\bj0  \otimes j a\bj1  &=  a\bj0  \otimes  a\bj1 j     ,
\end{align}
$j \in J$, $a,b  \in H$; \eqref{eq:smash2-coprod} and \eqref{eq:smash3-coprod} say that $H$ is a $J$-comodule algebra.
Then $J\ltimes H := J\otimes H$ with the tensor product structure
as an algebra and with the smash coproduct
\eqref{eq:smash coproduct-right} for the coalgebra structure is a Hopf algebra;
see e.g. \cite[1.1.4]{AN}\footnote{In loc. cit a left version is presented, with a different notation. The proof is equally straightforward.}.
We shall say that $H$ is a \emph{$J$-comodule Hopf algebra}.

\subsection{}\label{subsec:duals-co}

Let  $H$ and $U$ be  Hopf algebras, provided with a non-degenerate Hopf pairing
$(\, \vert \,): H \ot U \to \ku$.
\begin{lemma}\label{pairinglemma:co}
        Assume that $H$ and $U$ are $J$-comodules and that the pairing is compatible
        with $J$-coaction involving  the antipode of $J$, i.e.
        \begin{align} \label{eq:duality-modL-co}
        (a\bj0  \vert u)a\bj1 &= (a  \vert u\bj0)\Ss(u\bj1), & a \in H,\,  u \in U.
        \end{align}
        Then $H$ is a $J$-comodule Hopf algebra if and only if $U$ is so.
\end{lemma}

\pf
Let $u,v  \in U$, $a, b \in H$.  We compute
\begin{align*}
( (ab)\bj0 \vert u) (ab)\bj1 &= ( ab \vert u\bj0) \Ss(u\bj1)= (a\vert (u\bj{0})\_2) (b\vert (u\bj{0})\_1)  \Ss(u\bj1);
\\
(a\bj{0} b\bj{0} \vert u) a\bj1b\bj1
 &= (a\bj{0}\vert u\_2) (b\bj{0}\vert u\_1) a\bj1b\bj1
= (a\vert (u\_2)\bj{0}) (b\vert (u\_1)\bj{0}) \Ss((u\_2)\bj{1}) \Ss((u\_1)\bj{1})
\\
&=(a\vert (u\_2)\bj{0}) (b\vert (u\_1)\bj{0}) \Ss((u\_1)\bj{1}(u\_2)\bj{1}).
\end{align*}
Hence \eqref{eq:smash-coprod} holds for $U$  if and only if
 \eqref{eq:smash2-coprod} holds for $H$  and vice versa.
Similarly \eqref{eq:smash-counit} holds for $U$ if and only if \eqref{eq:smash-ep-coprod} holds for $H$ and vice versa.
Finally, \eqref{eq:smash3-coprod} holds for $H$ if and only if it holds for $U$:
\begin{gather*}
(a\bj0\vert u) j a\bj1 = (a\vert u\bj0) j \Ss (u\bj1) =  (a\vert u\bj0)  \Ss (u\bj1 \Sso(j));
\\
(a\bj0\vert u) a\bj1 j  = (a\vert u\bj0) \Ss (u\bj1) j =  (a\vert u\bj0)  \Ss (\Sso(j) u\bj1).
\end{gather*}
\epf

\subsection{}\label{subsec:co-smash-braided}

Let now $K$ be a Hopf algebra, $\toba$ a Hopf algebra  in $\kyd$ and also
a $J$-comodule coalgebra.
Extend the coaction of $J$ to  $H =  K \# \toba$ by
$\varrho(k \# b)  = k\# b\bj0 \ot b\bj1$, $b \in \toba$ and $k \in K$.
Then

\begin{itemize} [leftmargin=*]
        \item $H$ is a $J$-comodule coalgebra $\iff$ the coactions of $J$ and $K$ on $\toba$ commute, i.e.
        for all $b \in \toba$
\begin{align}\label{eq:coactions-commute}
(b\_0)\bj0 \ot b\_1 \ot (b\_0)\bj1 &= (b\bj0)\_0 \ot (b\bj0)\_1 \ot b\bj1 \in \toba \ot K \ot J.
\end{align}
        
        \item \eqref{eq:smash-ep-coprod} holds for $H \iff$ \eqref{eq:smash-ep-coprod} holds for $\toba$.
        Assume this is the case.
        
        \item \eqref{eq:smash2-coprod} holds for $H \iff$ \eqref{eq:smash2-coprod} holds for $\toba$ and
        the action of $k$ on $\toba$ is a morphism of $J$-comodules for all $k\in K$.
        
        \item \eqref{eq:smash3-coprod} holds for $H \iff$ \eqref{eq:smash3-coprod} holds for $\toba$.
\end{itemize}

In other words,  the coaction of $J$ on $H = K\# \toba$ satisfies \eqref{eq:smash-ep-coprod}, \eqref{eq:smash2-coprod}
and \eqref{eq:smash3-coprod}, if and only if so does the coaction of $J$ on $\toba$, and the
coaction of $J$ on $\toba$ commutes both with the action and the coaction of $K$. This can be phrased also as:
the homothety $\eta_\ell$ for $\ell\in J^*$ is a morphism of
Yetter-Drinfeld modules, i.e. $\eta_\ell \in \End^K_K \toba$.

\begin{definition}
        A  \emph{$J$-comodule braided Hopf algebra} is a  Hopf algebra $\toba$ in
        the braided category $\kyd$ that is also  a $J$-comodule coalgebra, that satisfies \eqref{eq:smash-ep-coprod}, \eqref{eq:smash2-coprod} and  \eqref{eq:smash3-coprod},  and such that the coaction of $J$ on $\toba$ commutes both with the action and the coaction of $K$. In such a case, the
        bosonization $H = K \# \toba$ is a $J$-comodule Hopf algebra
         and we can form the Hopf algebra $J \ltimes H = J \ltimes  (K \# \toba)$.
\end{definition}

\medbreak
As in subsection \ref{subsec:duals-co}, we consider the situation with
non-degenerate pairings; this time internal to the braided monoidal
category $\kyd$ instead of $\vect \ku$. Concretely,
let $\E$ be a Hopf algebra  in $\kyd$ provided with a non-degenerate bilinear
form  $(\, \vert \,): \toba \ot \E \to \ku$, and  extend it  by \eqref{eq:duality}
to a pairing $\toba\otimes\toba\otimes\E\otimes \E\to\ku$.
\begin{itemize} [leftmargin=*]\renewcommand{\labelitemi}{$\diamond$}
\item
The fact that the pairing is internal to the category $\kyd$
means that the
bilinear form $(\, \vert \,)$ is a morphism in the monoidal 
category $\kyd$, where $\ku$ is
endowed with the structure of a trivial Yetter-Drinfeld module.

\item We assume that for every $a, \widetilde a \in \toba$, 
$u, \widetilde u \in \E$,
the conditions \eqref{eq:dual1}, \eqref{eq:dual2} and 
\eqref{eq:dual3} of a Hopf pairing,
relating coproduct, product, unit and counit of $\toba$ and $\E$
hold.
\end{itemize}

Then we have in the braided category $\kyd$ exactly the same situation
we considered in Lemma \ref{pairinglemma:co} in the braided category
$\vect \ku$. The same calculations, this time in the
category $\kyd$, yield:

\begin{lemma}\label{pairinglemma:braided-co}
        Assume that both $\toba$ and $\E$ are $J$-comodules and that \eqref{eq:duality-modL-co} holds.
        Then $\toba$ is a $J$-comodule braided Hopf algebra, if and only if $\E$ is so. \qed
\end{lemma}

\subsection{}\label{subsec:smash-alggp}

Let $G$ be an affine algebraic group over $\ku$ and let $J = \ku[G]$ be
the algebra of functions on $G = \Alg(J, \ku)$.
Here we use the convention \eqref{eq:duality}, i.e.
\begin{align*}
\langle \gamma\eta, j\rangle &= \langle \gamma, j\_2\rangle  \langle \eta, j\_1\rangle, & \gamma, \eta& \in G.
\end{align*}
Thus, being a (right) $J$-comodule means being a rational (right) $G$-module: $m\cdot \gamma 
= m\bj0 \langle \gamma, m\bj1   \rangle$; which of course is equivalent to being rational left $G$-module.
So, in what follows we work with left rational modules.  
The conditions \eqref{eq:smash-coprod} and \eqref{eq:smash-counit}, respectively \eqref{eq:smash2-coprod} and \eqref{eq:smash-ep-coprod},
in the definition of  $J$-comodule Hopf algebra
just say that $G$ acts on $H$ by coalgebra, respectively algebra, automorphisms, while
\eqref{eq:smash3-coprod} is automatic by the commutativity of $\ku[G]$. 
We summarize our findings:

\begin{prop}
        Let $H$ be a Hopf algebra and let $G$ be an affine algebraic group acting rationally on $H$ by Hopf algebra maps.
        Then $H$ is a $\ku[G]$-comodule Hopf algebra and we can form $\ku[G] \ltimes H$. \qed
\end{prop}

\begin{remark} Since $J$ is commutative, $\GK (\ku[G] \ltimes H) = \dim G + \GK H$, 
see e.g. \cite[3.10]{KL}.
\end{remark}

\subsection{}\label{subsec:pointed-smash-co}
Let $K$ be a Hopf algebra and $V \in \kyd$, $\dim V < \infty$. Then  $\Aut_K^K (V)$ is an algebraic group, whose Lie algebra is $\End^K_K (V)$.
Every  morphism of algebraic groups $G \to \Aut_K^K(V)$ extends to
an action of $G$ on $T(V)$ by Hopf algebra automorphisms in $\kyd$; hence it descends to 
an action of $G$ on $\toba(V)$ by Hopf algebra automorphisms in $\kyd$.
It extends to an action of $G$ on $K \#\toba(V)$, trivially on $K$, 
giving rise to the Hopf algebra  $\ku[G] \ltimes (K \#\toba(V))$.
One may define analogous actions of these Hopf algebras from bosonizations of pre-Nichols or post-Nichols algebras,
or of deformations of Nichols algebras, provided that the ideals of defining relations are preserved
by the action of $G$.

\section{Hopf algebras arising from Nichols algebras of diagonal type}

\subsection{} \label{subsec:smash-nichols-diagonal}

Let $\theta \in \N$, $\I = \I_\theta=\{1,2,\ldots,\theta\}$. Denote
by $(\alpha_i)_{i\in \I }$ the canonical  basis of $\zt$.

Let $(V, c)$ be a braided vector space of diagonal type of dimension $\theta$;
let $(x_i)_{i\in \I }$ be a basis of $V$. Since $(V,c)$
is assumed to be of diagonal type,
there is a matrix $\bq = (q_{ij})_{i,j\in \I }\in (\ku^\times)^{\I \times \I}$
such that $c(x_i\otimes x_j)=q_{ij}x_j\otimes x_i$ for all $i,j\in \I$.
Then the tensor algebra $T(V)$ and and the Nichols algebra
$\toba(V)$ are $\zt$-graded (as braided Hopf algebras), by $\deg x_i = \alpha_i$, $i\in \I$.

Let $K$ be a Hopf algebra. To realize the braided
vector space $(V, c)$ as a Yetter-Drinfeld module over  $K$
we need some extra data.
\begin{itemize} [leftmargin=*]\renewcommand{\labelitemi}{$\diamond$}
\item
A pair
$(g, \chi) \in G(K)\times \Alg(K, \ku)$ is called 
a \emph{YD-pair} \cite{AAnGMV} if      
$\chi(a)\,g = \chi(a\_{2}) a\_{1}\, g\, \Ss(a\_{3})$ for all 
$a\in K$. This implies $g\in Z(G(K))$.
\item
Then   $\ku_g^{\chi} := \ku$ with coaction given by $g$ and action given by $\chi$ is a simple object in $\ydk$.
\end{itemize}
A \emph{principal realization} of the braided vector space 
$(V,c)$ over the Hopf algebra $K$ is a family
$\left((g_i, \chi_i)\right)_{i \in \I}$ of YD-pairs such that
\begin{align}\label{eq:realization-ppal-diag}
\chi_j(g_i) &= q_{ij},& \text{for all }&  i,j\in \I.
\end{align}

A principal realization allows us to see  braided vector space
as a Yetter-Drinfeld module, $V\in \ydk$, by 
declaring $x_i \in V_{g_i}^{\chi_i}$, $i\in \I$.
Let $ d_{g}^{\chi} = \dim V_{g}^{\chi} = \vert \{i\in \I: (g_i, \chi_i) = (g, \chi)\}\vert$. Then
\begin{align*}
\bd_V=\End_K^K(V) & \simeq \bigoplus_{g\in \Gamma, \chi \in \widehat\Gamma} \gl(d_{g}^{\chi}, \ku).
\end{align*}

Despite the notation, the Lie algebra $\bd_V$ depends on the way the
braided vector space $V$ is realized as a
$K$-Yetter-Drinfeld module and not merely on the braided vector space
$V$ itself.

\medbreak
For $h = (h_i)_{i\in \I_{\theta}}\in \ku^{\theta}$ we denote by
$D_{h} \in \End(V)$  the map defined by
$D_{h}(x_i)= h_ix_i$, $i \in \I_{\theta}$.
By abuse of notation, we denote by $D_{h}$ the corresponding derivation of
$T(V ) \# \ku\Gamma$ or $\toba(V) \# \ku\Gamma$. Let
\begin{align*}
\tg_V & = \{D_{h}: h\in \ku^{\theta}\} \subseteq \bd_V.
\end{align*}
The abelian Lie algebra $\tg_V$ depends only on $(V,c)$.
If $(g_i, \chi_i) = (g_j, \chi_j)$ implies $i = j$, then $\bd_V = \tg_V$.

\begin{remark}\label{rem:t-homogeneous}
The action of the Lie algebra $\tg_V$ preserves the $\zt$-grading. Indeed, let  $h \in \ku^{\theta}$ and let $\alpha \mapsto h_{\alpha}$
be the unique group homomorphism  $\zt \to \ku$ such that $h_{\alpha_i} = h_{i}$, $i\in \I$.
Then $D_{h}$ acts by $h_{\beta}$ in the homogeneous component $T(V)_\beta$ for all $\beta \in \zt$.
Hence every Hopf ideal $\cI$ of $T(V)$ generated by $\zt$-homogeneous elements is stable under $\tg_V$
and $\tg_V$ acts by derivations and coderivations on $\T(V)/ \cI$.
\end{remark}

\begin{remark}
In fact, the $\zt$-grading is tantamount to a comodule structure over the group algebra $\ku \zt$,
which is the algebra of functions on the algebraic torus $\Tb_V$; $\tg_V$ is its Lie algebra, and the action of
$\tg_V$ is the derivation of the natural action of $\Tb_V$.
\end{remark}

\subsection{} \label{subsec:smash-double-diagonal} From now on,
we assume that  $\car \ku =0$. We keep the notation above and assume that
$\dim \toba(V) < \infty$. The classification of the finite-dimensional Nichols algebras of diagonal type 
was given in \cite{H-classif}.
An efficient set of defining relations of $\cB(V)$, i.e. generators of the ideal $\cJ_{\bq}$,
was provided in \cite{Ang-crelle}. Besides $\toba(V)$, there are
two other  Hopf algebras in $\ydk$ that are expected to play a  role in
 representation theory:

\begin{enumerate}[leftmargin=*,label=\rm{(\alph*)}]
\item\label{item:prenichols}
\cite{Ang-crelle,A-pre-Nichols} The \emph{distinguished pre-Nichols algebra} of $(V, c)$ is the quotient $\dpnq := T(V)/\cI_{\bq}$ by a suitable ideal $\cI_{\bq}$. Thus, 
there are projections $T(V) \twoheadrightarrow \dpnq \twoheadrightarrow \toba(V)$.

\medbreak\item\label{item:postinichols} \cite{AAR} The \emph{Lusztig algebra} of $(V, c)$ is
 the graded dual $\luq$ of $\dpnq$.  
        \end{enumerate}

\begin{prop}\label{lema:unrolled-prenichols}
 Let $K$ be a Hopf algebra provided with a principal realization of $(V,c)$ and let $L= U(\tg_V)$.
 Then $\dpnq$ and $\luq$ are $L$-module braided Hopf algebras in $\ydk$ and we can form the unrolled bosonizations $(\dpnq\#K) \rtimes L$ and
 $(\luq\#K) \rtimes L$.
 \end{prop}

 \pf
 The claim for $\dpnq$ follows from Remark \ref{rem:t-homogeneous} and  implies the one for $\luq$ by
 Lemma \ref{pairinglemma:braided}.
 \epf

\begin{example}
	If $\theta = 1$ and $\bq$ is a root of 1 of even order, then we  recover the construction in \cite{GPT, CGP}.
\end{example}

 \bigbreak
 \subsection{}\label{subsec:unrolled-post-nichols}
 
 Let $(V,c)$ be of diagonal type with $\dim \toba(V) < \infty$. 
 Fix a principal realization over the group algebra $\ku \Gamma$, where $\Gamma$ is abelian.
 Then  each of the Hopf algebras $\toba(V)$, $\dpnq$ and $\luq$ in $\ydG$ gives rise to Hopf algebras
 $\ug(V)$, $U(V)$, $\Uc(V)$ respectively; they are suitable Drinfeld doubles of the bosonizations
 $\toba(V)\# \ku\Gamma$, $\dpnq\# \ku\Gamma$ and $\luq\# \ku\Gamma$. See \cite{H-lusztig iso,A-pre-Nichols,AAR}.
 If $\bq$ is symmetric, then we may divide that Drinfeld double by a central Hopf subalgebra. 
 If furthermore $\bq$ is of Cartan type, then we recover the small, the De Concini-Procesi 
 and the Lusztig divided power quantum group, respectively. Then we may define unrolled quantum groups
\begin{align}
\uq &\rtimes U(\tg_V),& U(v) &\rtimes U(\tg_V),& \Uc(V) &\rtimes U(\tg_V).
\end{align}
Indeed, the Lie algebra $\tg_{V \oplus W}$ acts on $T(V \oplus W) \# \ku\Gamma$, but 
if  $\zeta \in \ku^{2\theta}$, then
$D_{\zeta}$ preserves the relations of the quantum double
if and only if $\zeta$ belongs to the image of the map $\tg_{V} \to \tg_{V \oplus W}$,
$\xi \mapsto (\xi, -\xi)$.

\end{document}